\documentclass[10pt]{amsart} 

\usepackage{amsfonts,amsmath,amsthm, graphics}

\theoremstyle{plain}
\newtheorem{theorem}{Theorem}

\newtheorem{definition}{Definition}

\title[R.F. Shamoyan, M.G. Bashmakova]{On embedding theorems in analytic and harmonic function spaces of several variables \\in some domains in  $\mathbb{C}^n$}

\author{R.F. Shamoyan, M.G. Bashmakova}

\address{Department of Mathematics, Bryansk State Technical University, Bryansk, 241050, Russia}

\email{\rm rshamoyan@gmail.com, 
mariya-bashmakova@yandex.ru}

\begin{document}

\begin{abstract}
In this  survey we collect some recent advances concerning embedding theorems in analytic and harmonic function spaces of several variables in various domains in $\mathbb{C}^n.$

Some sharp embedding results presented in  this survey paper extend sharp embeddings in analytic function spaces obtained in higher dimension previously by J. Ortega, J. Fabrega, J. Cima, D. Luecking, M. Abate and many others in recent decades in the unit ball and bounded strongly pseudoconvex domains with the smooth boundary for onefunctional case to similar type analytic multyfunctional Bergman type function spaces. Same type multifunctional  sharp extensions of some recent embedding results of one functional analytic function spaces of B. Sehba and D. Bekolle will be  provided by us in this survey for unbounded domains namely for tubular domains over symmetric cones.

In this paper we collect many new sharp embeddings in higher dimension obtained recently by first author though some one dimensional new sharp embeddings will be also mentioned shortly at the end of paper.

Our results can be seen as direct sharp extensions of various sharp embedding theorems in analytic spaces in the unit disk obtained many years ago by L. Carleson, P. Duren, I Verbitsky, D. Luecking and many others.

Nice properties of so called $r$-lattices in various domains in $C^n$  are very important in all proofs of theorems of this paper.

\end{abstract}

\maketitle

\section{Introduction}
The study of embedding theorems in analytic function spaces in one variable for Hardy spaces started from a celebrated paper of L.~Carleson \cite{k19}, then similar type sharp embedding  resutls for other function spaces were obtained in papers of B.~Oleinik, B.~Pavlov \cite{k6}, \cite{k7}, N.~Shirokov \cite{shir}, P.~Duren \cite{k18}, F.~Shamoyan \cite{d8}-\cite{d7}, W.~Hastings \cite{k5} and I.~Verbitsky \cite{v1} and many others in the unit disk or in the unit polydisk. Later some very interesting  work in this direction were presented in several papers by W.~Cohn, I.~Verbitsky, P.~Ahern \cite{d49}, J.~Cima \cite{d5} and D.~Luecking \cite{k17}, \cite{d17} and their coauthors and others in the context of the unit ball. L.~Hormander \cite{k20} completely extended Carlesons
sharp result for Hardy class to very difficult domains namely to convex domains of finite type. We indicated partially  these papers in the list of references of this survey. Referying interested  readers also to the list of referenes of all these papers. We mention separately that
Igor Verbitsky in 1987  completely solved such type embedding problems providing a series of sharp results in context of so called analytic mixed norm spaces in the unit disk using very delicate techniques of vector valued inequalities and dyadic partitions of the unit disk and delicate estimateds obtained earlier by E. Stein and coauthors for tent spaces in several variables. We refer the reader to this paper \cite {v1} and to other many interesting other papers concerning embeddings  which I. Verbitsky mentioned in this work.

In function spaces in $\mathbb{R}^n$ such type numerous embedding results complete analogues of mentioned theorems were also considered by many authors in  recent several decades. We don't discuss  these issues however here in this paper.

The first author in recent decades obtained a large amount of new sharp embedding theorems in various function spaces for functions of several variables in various domains and we plan to collect these results in this expository paper. It is well known that many sharp interesting embedding theorems in one and several complex variables in various complicated domains have many nice applications in complex function theory.
Starting from the study of various aspects of Toeplits operators and multipliers (pointiwise or coefficient), corona theorem and in various
other issues in complex function theory and operator theory.

First we very shortly review some recent results concerning embeddings in analytic area Nevanlinna type spaces of several complex variables which is a new topic.
Namely we consider here the very simple domain the unit polydisk, though more difficult domains as tube or bounded strongly pseudoconvex domains may be also considered here.

Here in this completely new research area new nice open questions may be posed related with the embeddings in area Nevanlinna spaces of several variables. We leave this task to interested readers.

Then we pass to various new interesting multyfunctional sharp embedding results in various domains in $\mathbb{C}^n$ in various interesting  function spaces of several variables.
As far as we know this topic is also new.

Previously such type sharp embeddings almost new results were known only for one function spaces (onefunctional spaces). Though many sharp results were given here in difficult domains and in various function spaces of several complex variables. Further we consider mixed norm and analytic Herz spaces.
 Here to be more precise various new sharp embeddings in so called new mixed norm spaces will be discussed  and in so called Herz type spaces of analytic function spaces of several variables.  But with only one function within the inequality of the embedding though in various difficult domains.

Further we discuss and provide shortly certain new sharp embeddings in another new research area namely we consider  new embedding theorems which are sharp and which are related with the so called TRACE function in tubular domains over symmetric cones and in bounded strongly pseudoconvex domains.
Such type sharp embedding theorems in comlpex function space theory of several variables are completely new and much work should be done here.
We consider in this work sharp embedding theorems in various analytic function spaces in various domains, namely unit polydisk, unit ball, tubular domains over symmetric cones and bounded strongly pseudoconvex domains with the smooth boundary.

All these interesting new results can be seen in many recent papers of the first author and we collected all of them in this expository research paper.

We collected below at the end of paper a very short review of recent research papers on embedding theorems in analytic function spaces of several variables. Since sizes of this paper will not allow to present all results concerning this topic we wanted.

In many theorems below we define and work with Herz spaces of analytic and harmonic functions of several variables. We define two such function spaces.of several variables. Note that the role of various properties of so called $r$-lattices is very important in many proofs of our sharp embedding  theorems.

\section{Embedding theorems in some domains in $C^n$}

Let $U^n$ be the unit polydisk in $\mathbb{C}^n$, $U^n=\{z\in\mathbb{C}^n : |z_j|<1, j=1,...n\}$. Let $H(U^n)$ be the space of all analytic functions in $U^n.$ Let further $\mathbf{B}_n$ be the unit ball in $\mathbb{C}^n, \mathbf{B}_n=\{z\in\mathbb{C}^n:|z|<1\}, S_n=\{z: |z|=1\}$ be unit sphere. Let $H(\mathbf{B}_n)$ be the space of all
analytic functions in $\mathbf{B}_n.$ We define by $d\nu_\alpha$ as usual weighted Lebegues measure depending on domain in $\mathbb{C}^n.$

In this paper everywhere below we denote by $\mu$ positive Borel measure on $D$ domain where $D$ is a general domain in $C^n.$

We define in many places below   $\Delta_k$ as $B(a_k, r)$ where $B$ is a usual standard Bergman ball in various $D$ domains in higher dimension and $a_k$ is a $r$-lattice in such domains in $C^n$.

Let further
$$N_{\omega}^{p,q}(U^n)=\{ f\in H(U^n) : \int\limits_{T^n}\big(\int\limits_{I^n}(ln^+|f(\tau_1\xi_1,...\tau_n \xi_n)| )^p\prod\limits_{k=1}^n \omega(1-\tau_k)d\tau_k\big)^{\frac{q}{p}}d\xi_1...d\xi_n<+\infty \},$$
$$\widetilde{N}_{\omega}^{p,q}(U^n)=\{ f\in H(U^n) : \int\limits_{T^n}\big( \int\limits_{I^n}(ln^+|f(\tau_1\xi_1,...\tau_n \xi_n)| )^p d\xi_1...d\xi_n \big)^\frac{q}{p} \prod\limits_{k=1}^n \omega(1-\tau_k)d\tau_k <+\infty  \},$$
where as usual
$$T^n=\{z\in\mathbb{C}^n: |z_j|=1, j=1,...,n \}, I^n=[0,1]^n, $$ and $0<p,q<\infty, $ where $u^+=\max(u,0).$ Our weight can be also radial $(1-r)^\alpha.$ In this case we also define $ N_{\overrightarrow{\omega}}^{p,q}, N_{\alpha}^{p,q}=N_\alpha^p,$ if $\omega(r)=r^\alpha, \alpha>-1, p=q.$\\

Let $\Lambda$ be a bounded (or unbounded) domain with $\mathbb{C}^2,$  let $H(\Lambda)$ be the space
of all analytic spaces in $\Lambda$. We denote by $dV_\alpha$ the weighted Lebegues measure on $\Lambda$. Let
$\partial \Lambda$ be boundary of $\Lambda$. Let further
$$dV_\alpha(z)=( \textit{dist}(z,\partial \Lambda))^{\alpha}dV(z), z\in \Lambda, \alpha>-1,$$
where $dV$ is a normalized Lebegues measure on $\Lambda$.
Further we considered also area Nevanlinna spaces on $\Lambda$
$$N_{\alpha}^p(\Lambda)=\{f\in H(\Lambda): \|f\|_{N_\alpha^p(\Lambda)}^p=\int\limits_{\Lambda}(\log^+|f(z)|)^p dV_\alpha(z)<\infty \}, 0<p\leq \infty,$$
and these are Banach spaces for $p \geq 1$ and complete metric spaces for other values of $p.$

For many interesting properties of a special $\omega$ function
which appear below in two interesting sharp theorems below.
We refer the reader to \cite{1}  namely for weighted $\omega$ function spaces from special  $S$ function class. We refer the reader to the list of references of \cite{1} for more info about that class  of interesting special weights also.
Various applications of $\omega$ function class of these special weights  also can be seen in mentioned papers.
In the following sharp embedding theorem in the unit polydisk we provide new interesting sharp results on the action of so called differentiation $D$ operator in area Nevanlinna spaces of several variables.
In more general domains which we defined above these results may also be valid. We pose this interesting question as a problem for interested readers.
Note that in \cite{1} we can see a long interesting history of this nice problem also on the action of differentiation operator in so called area Nevanlinna spaces of one variable.

Though we formulate all our sharp embedding results, theorems in this paper in one domain in all theorems below our proofs are not very complicated and may be probably passed to other domains by similar arguments. We leave this task to interested readers.

We denote by $D$ stadard differentiation operator (see \cite{1}).

\begin{theorem} (see \cite{1}). 
Let $0<p<+\infty, \int\limits_0^{1}\omega_j(t)dt<+\infty,j=1,2...,n.$
Then
$$\int\limits_{I^n}(\int\limits_{T^n}ln^{+}|Df(\tau_1\xi_1,...,\tau_n\xi_n)|d\xi_1...d\xi_n)^p\prod\limits_{j=1}^n\omega_j(1-\tau_i)d\tau_1...d\tau_n\leq$$
$$ \leq C\int\limits_{I^n}(\int\limits_{T^n}ln^+|f(\tau_1\xi_1,...\tau_n\xi_n)|d\xi_1...d\xi_n)^p\prod\limits_{j=1}^n\omega_j(1-\tau_i)d\tau_1...d\tau_n,$$
if and only if
$$ \int\limits_0^1 \omega_j(t)(ln\frac{1}{t})^pdt<+\infty, j=1,2,...n.$$
\end{theorem}

\begin{theorem} (see \cite{1}). 
Let $s\geq 1, s\geq \max (q,p), \omega=\prod\limits_{j=1}^n w_j.$ Let
$$\frac{2}{s}-\frac{1}{p}>0, \widetilde{\omega}_j(1-|z_j|)=\omega_j(1-|z_j|)^{\frac{q}{s}}(1-|z_j|)^{\frac{2q}{s}-\frac{q}{p}-1}.$$
Then $Df$ is acting from $N_{\widetilde{\omega}}^{p,q}, \widetilde{N}_{\widetilde{\omega}}^{p,q}$ to $N_{\omega}^{s,s}$ if and only if
$$\int\limits_{0}^{1}\omega_j(1-\tau)ln\big(\frac{1}{1-\tau}\big)^s d\tau <+\infty, j=1,2,...n.$$
\end{theorem}

Below everywere we denote by $T_\Omega$ tubular domains over summetric cones, by $H(T_\Omega)$ we denote spaces of all analytic functions in $T_\Omega$ domains and by $A^p_s$, $A^{pq}_s$ we denote the well known analytic Bergman and so called mixed norm spaces of analytic functions on these unbounded tubular domains over symmetric cones. We refer for precise definitions of these objects to papers of B. Sehba, D. Bekolle  (see, for example, \cite{k4}) and his various coauthors in the list of references of this expository paper.

In the following five theorems we collect sharp new embedding theorems for new multifunctional analytic Bergman and Herz spaces in rather general and complicated tubular domains over symmetric cones in function spaces of several complex variables.
The interested reader should look at \cite{2} for many simple (or not) and many  very basic definitions of unbounded complicated tubular domains over symmetric cones and many basic definitions of analytuc function spaces on them closely related to this interesting new valuable  research area.
Analytic function spaces of several variables (for example, classical Bergman spaces) on these domains were studied in recent papers of A.~Bonami, D.~Bekolle, B.~Sehba and many others (see, for example, \cite{k4}).
We define multifunctional extensions below of these classes and show sharp new embeddings.

Let $A$ be a Banach space of holomorphic functions
a domain $T_{\Omega}\in \mathbb{C}^n;$
 given $p\geq 1$, a finite positive Borel measure $\mu$ on $T_{\Omega}$ is a Carleson measure
of $A$ (for p) if there is a continuous inclusion $A\rightarrow L^p(\mu)$, that is there exists a constant $C > 0$ such that

$$\forall f\in A  \int\limits_{T_\Omega}|f|^p d\mu\leq C \|f\|_A^p.$$

 We define new analytic general multifunctional Herz (for $p = q$ Bergman) spaces
in tubular domains over symmetric cones with smooth boundary as spaces of $(f_k)_{k=1}^m$ functions analytic in $T_\Omega$, so that
$$||(f_1,...,f_m)||^{q}_{A(p,q,m,d\mu)}=\sum\limits_{k=1}^{\infty}\left( \int\limits_{B(a_k,r)}\prod\limits_{i=1}^m |f_i(z)|^p d\mu(z)\right)^{\frac{q}{p}}< \infty, 0<q,p<\infty,$$
$\mu$ is a positive Borel measure on $T_\Omega$.

We can very easily define very similarly many such type Herz type analytic function spaces in many domains in $C^n$ as far as $r$-lattices there exits. In the papers of first author with variuos  coauthors in recent two decades many sharp embeddings were proved for such spaces in various domains though using similar approaches and theniques and we list them in various theorems below.

Another type of analytic Herz type function spaces in various domains in $C^n$ can be defined in higher dimension also if we simply replace summation in definition above by integration, such type analytic function  spaces also appear often below.

Herz spaces can be defined if we replace summation by integration, all our results may be valid there.

We refer for $\delta(a_k)$ to \cite{2}. It is a Lebeques measure of  $B(a_k,r).$

\begin{theorem} (see \cite{2}). 
Let $0<p,q<\infty, 0<s\leq q, \beta_j>-1$ for $j=1,...m,$ and let $\mu$ be a
positive Borel measure on $T_\Omega$. Then the following conditions are equivalent:

1. If $f_i\in A^s_{\beta_i}, i=1,...m$ then
$$||(f_1,...f_m)||_{A(p,q,m,d\mu)}\leq C\prod\limits_{i=1}^m ||f_i||_{A^s_{\beta_i}} .$$

2. The measure $\mu$ satisfies a Carleson type condition:
$$ \mu(\Delta_k)\leq C\delta(a_k)^{\sum\limits_{i=1}^m\frac{p(2n/r+\beta_j)}{s}}, k\geq 1,$$
for any $\{a_k\}$ $r$-lattice in $T_\Omega$.
\end{theorem}

\begin{theorem} (see \cite{2}). 
Let $0<p,q<\infty, 0<s\leq q, \beta_j>-1,$ for $j=1,...m, 0<t_i\leq s,$ for $i=1,...m, \sum\limits_{i=1}^n\frac{1}{t_i}=\frac{m}{s}, m\in \mathbb {N}$,  and let $\mu$ be a positive Borel measure on $T_\Omega$. Then the following conditions are equivalent:

1. If $f_i\in A^{s,t_i}_{\beta_i}, i=1,...m$ then
$$||(f_1,...f_m)||_{A(p,q,m,d\mu)}\leq C\prod\limits_{i=1}^m ||f_i||_{A^{s,t_i}_{\beta_i}} .$$

2. The measure $\mu$ satisfies a Carleson type condition:
$$ \mu(\Delta_k)\leq C\delta(a_k)^{\sum\limits_{i=1}^m\frac{p(2n/r+\beta_j)}{s}}, k\geq 1,$$
for any $\{a_k\} $ $r$-lattice in $T_\Omega$.
\end{theorem}

We denote by $B(w,r)$ the Bergman ball of center $w$  and radius  $r$. We define as usual by $d\nu_\alpha$ weighted Lebegues measure in tube domain.

We refer the reader for various objects in theorems below to related references in those theorems and we refer  to \cite{2} for definition of $\delta (z)$ function.

\begin{theorem} (see \cite{2}). 
Let $0<p,q<\infty, 0<\sigma_i\leq q, -1<\alpha_i<\infty$ for $i=1,...,m.$  Let $\mu$ be a
positive Borel measure on $T_\Omega$. Then the following two conditions are equivalent:

1. $||(f_1,...f_m)||^q_{A(p,q,m,d\mu)}\leq C\prod\limits_{i=1}^m\int\limits_{T_\Omega}\left(\int\limits_{B(w,r)}|f_i(z)|^{\sigma_i}\delta^{\alpha_i}(z)d\nu(z)\right)^{\frac{q}{\sigma_i}}dw$

2. The measure $\mu$ satisfies a Carleson type condition:
$$ \mu(\Delta_k)\leq C\delta(a_k)^{m(2n/r)\frac{p}{q}+\sum\limits_{i=1}^m\frac{p(2n/r+\alpha_i)}{\sigma_i}}, k\geq 1$$
for any $\{a_k\}$ $ r$-lattice in $T_\Omega$.
\end{theorem}

\begin{theorem} (see \cite{2}). 
Let $\mu$ be a positive Borel measure on $T_\Omega$. Let $0 < p_i,q_i<\infty, i = 1,...,m$, let $\sum\limits_{i=1}^m\frac{1}{q_i}=1, \alpha>-1.$ Then the following conditions are equivalent:

1. If $f_i\in H(T_\Omega), i=1,...m,$ then $$ \int\limits_{T_\Omega}\prod\limits_{i=1}^m|f_i(z)|^{p_i}d\mu(z)\leq c\prod\limits_{i=1}^m \left( \sum\limits_{i=1}^{\infty}\left(\int\limits_{\Delta_k^*}|f_i(z)|^{p_i}\delta^\alpha(z)d\nu(z)\right)^{q_i} \right)^{\frac{1}{q_i}}.$$

2. The measure $\mu$ satisfies a Carleson type condition $$\mu(\Delta_k)\leq C|\Delta_k|^{m(1+\frac{\alpha}{2n/r})},$$ for $k\geq 1.$
\end{theorem}

\begin{theorem} (see \cite{2}). 
Let $0<p_i,\sigma_i<\infty, -1<\alpha_i<\infty, i=1,...,m.$ Let $\mu$ be a positive Borel
measure on $T_\Omega$. Then the following conditions are equivalent:

1. If $f_i\in H(T_\Omega), i=1,...m,$ then $$ \int\limits_{T_\Omega}\prod\limits_{i=1}^m|f_i(z)|^{p_i}d\mu(z)\leq c\prod\limits_{i=1}^m \int\limits_{T_\Omega}\left( \int\limits_{B(w,r)}|f_i(z)|^{\sigma_i}d\nu_{\alpha_i}(z)\right)^{\frac{p_i}{\sigma_i}}d\nu(w). $$

2. The positive Borel measure $\mu$ satisfies the following Carleson type condition: $\mu(\Delta_k)\leq c_1\delta(a_k)^\tau,$ where \\ $\tau=m(2n/r)+\sum\limits_{i=1}^{m}\frac{p_i(2n/r)+\alpha_i}{\sigma_i}, k\geq 1, $ for some constants $c$
and $c_1$, and for any $\{a_k\}$ $r$-lattice in $T_\Omega$.
\end{theorem}

Results which we presented in theorems 3-7 may be valid with same or very similar proofs in various other domains also. We leave this to interested readers.

In the following new interesting sharp several theorems we collected various new  sharp embedding theorems in harmonic function spaces of several complex variables obtained recently by the first author.in ccoperation with Milos Arsenovic.
These series of new sharp interesting assertions provided by us below 
in spaces of harmonic functions of several variables can be considered as complete analogues of our previous sharp embedding theorems for analytic functions in other domains in tubular domains over symmetric cones  which we provided above and bounded strongly pseudoconvex domains which will be formulated by us below.
Spaces of harmonic functions which we consider below are completely  new and these harmonic multyfunctional function spaces can be considered as direct extentions of well studied analytic Bergman function spaces of one and several complex  variables.
We refere the reader for various simple well known and basic definitions of harmonic function spaces theory to \cite{3}  and also to various references there also. Various interesting aspects of 
harmonic function spaces of several variables studied in recent years by various authors.
We refer the reader for example to recent  papers of A. Alexandrov and his coauthors and to papers which can be seen in the list of references of mentioned papers.
This research area is relativly new and many open interesting problems can be seen here.

We set $\mathbb{H}=\{(x,t): x\in \mathbb{R}^n, t>0\}\in \mathbb{R}^{n+1}.$ For $z=(x,t)\in \mathbb{H}$ we set $\overrightarrow{z}=(x,-t).$

The space of all harmonic functions in a domain $\Omega$ is denoted by $h(\Omega)$. Weighted
harmonic Bergman spaces on $\mathbb{H}$ are defined, for $0 < p < \infty$ and $ \lambda>-1,$ by
$$ A(p,\lambda)=A(p,\lambda)(\mathbb{H})= \{ f\in h(\mathbb{H}):||f||_{A(p,\lambda)}=\left( \int\limits_{\mathbb{H}} |f(z)|^{p}d\mu_{\lambda}(z)\right)^{\frac{1}{p}}<\infty \}.$$

For $f\in h(\mathbb{H})$ and $t>0$ we set $M_p(f,t)=\|f(\cdot,p)\|_{L_p(\mathbb{R}^n)}, 0<p<\infty,$ with the usual
convention in the case $p=0$. We use harmonic mixed norm spaces

$$B_{\alpha}^{p,q}=\bigg\{ f\in h(\mathbb{H}):\|f\|^{q}_{B^{\alpha}_{p,q}} =\int\limits_0^{\infty}M_p^q(f,t)t^{\alpha q-1}dt<\infty \bigg\},$$ where $0\leq p <\infty, 0<q<\infty,$ and harmonic Triebel-Lizorkin spaces
$$F_{\alpha}^{p,q}=\bigg \{f\in h(\mathbb{H}):\|f\|^{p}_{F_{\alpha}^{p,q}}=\int\limits_\mathbb{R} \bigg( \int\limits_0^{\infty}|f(x,t)|^q t^{\alpha q-1}dt \bigg)^\frac{p}{q} dx<\infty\bigg\},$$
where again $0\leq p \leq \infty, 0<q\leq\infty,$  the case $q=\infty$ is covered by the usual convention.

By $X-L^p$ Carleson measure (or simply $X$ Carleson measure when $p$ is clear from the
context) of a (quasi)-normed subspace $X$ of $h(\Omega)$ we understand those positive Borel
measures $\mu$ on $\Omega$ such that
$$ \bigg (\int\limits_{\Omega}|f|^{p}d\mu \bigg)^{\frac{1}{p}}<C\|f\|_{X},   f\in X. $$

Typical cases are $\Omega=\mathbb{B}=\{x\in \mathbb{R}^n :|x|<1 \}$  and $\Omega=\mathbb{H}.$

\begin{theorem} (see \cite{3}). 
Let $\mu$ a positive Boreal measure on $\mathbb{H}$. Assume $0<p_i, q_i<\infty, i=1,..,m,$ satisfy $\sum\limits_{i=1}^m\frac{1}{q_i}=1,$ and let $\alpha>1$.Then the following two conditions on the measure $\mu$ are equivalent:

1. If $f_i, i=1,...,m,$ are functions harmonic in $\mathbb{H}$, then we have
$$\int\limits_{\mathbb{H}}\prod\limits_{i=1}^m |f_i(z)|^{p_i}d\mu(z)\leq C\prod\limits_{i=1}^m \bigg [ \sum\limits_{k=1}^{\infty}\bigg( \int\limits_{\Delta^*_k}|f_i(z)|^{p_i}t^{\alpha}dz\bigg)^{q_i}\bigg]^{\frac{1}{q_i}}.$$

2. The measure $\mu$ satisfies a Carleson type condition:
$$\mu(\Delta_k)\leq C |\Delta_k|^{m(1+\frac{\alpha}{n+1})}, k\geq 1.$$
\end{theorem}

\begin{definition}
Let $0<q,p<\infty$ and let $\mu$ be a positive Borel measure on $\mathbb{H}$. The space $A(p,q,md\mu)$ is the space of all $(f_1,...,f_m)$ where $f_i\in h(\mathbb{H}), i=1,..,m,$ such that
$$ \|(f_1,...,f_m)\|^q_{A(p,q,md\mu)}=\sum\limits_{k=1}^{\infty}\bigg( \int\limits_{\Delta_k}\prod\limits_{i=1}^m|f_i(z)|^pd\mu(z)\bigg)^{\frac{q}{p}}<\infty.$$
\end{definition}

If $d\mu=dm_s,$ then we write $A(p,q,m,s)$ for the corresponding space, if $m = 1$ we write $A(p,q,d\mu)$ or $A(p,q,s)$ if $d\mu=dm_s.$

We refer the reader to \cite{3} for definition of $\eta_k$. It is a Lebeques measure of  $B(a_k,r)$.

\begin{theorem} (see \cite{3}). 
Let $0<p,q<\infty, 0<s\leq q, \beta_i\geq -1,$ for  $i=1,...,m$ and let $\mu$ be a
positive Borel measure on $\mathbb{H}$. Then the following conditions are equivalent:

1. If $f_i\in A(s,\beta_1), i=1,...,m,$ then
$$\|(f_1,...,f_m)\|_{A(p,q,m,d\mu)}\leq C\prod\limits_{i=1}^m\|f\|_{A(s,\beta_i)}.$$

2. The measure $\mu$ satisfies a Carleson type condition:
$$\mu(\Delta_k)\leq C\eta_k^{\sum\limits_{i=1}^m\frac{p(n+1+\beta_i)}{s}}, k\geq 1.$$
\end{theorem}

\begin{theorem} (see \cite{3}). 
Let $0<p,q<\infty, 0<s\leq q, \beta_i\geq -1,$ for  $i=1,...,m, 0<t_i<s$ for $i=1,...,m$ and let $\mu$ be a positive Borel measure on $\mathbb{H}$. Then the following conditions are equivalent:

1. If $f_i\in B^{s,t_i}_{(\beta_i+1)/s}, i=1,...,m,$ then
$$ \|(f_1,...,f_m)\|_{A(p,q,m,d\mu)}\leq C\prod\limits_{i=1}^m \|f\|_{B^{s,t_i}_{(\beta_i+1)/s}}.$$

2. If $f_i\in F^{s,t_i}_{(\beta_i+1)/s}, i=1,...,m,$ then
$$ \|(f_1,...,f_m)\|_{A(p,q,m,d\mu)}\leq C\prod\limits_{i=1}^m \|f\|_{F^{s,t_i}_{(\beta_i+1)/s}}.$$

3. The measure $\mu$ satisfies a Carleson type condition:
$$\mu(\Delta_k)\leq C\eta_k^{\sum\limits_{i=1}^m\frac{p(n+1+\beta_i)}{s}}, k\geq 1.$$
\end{theorem}

We refer for $Q_w$ to \cite{3}.

\begin{theorem} (see \cite{3}). 
Let $0<p,q<\infty $ and $0<\sigma_i<q, -1<\alpha_i<\infty$ for $i=1,...,m.$ Let $\mu$ be
a positive Borel measure on $\mathbb{H}$. Then the following two conditions are equivalent:

1. For any $m$-tuple $(f_1,..., f_m)$ of harmonic functions on $\mathbb{H}$ we have
$$\|(f_1,...,f_m)\|^q_{A(p,q,m,d\mu)}\leq C\prod\limits_{i=1}^m \int\limits_{\mathbb{H}}\bigg ( \int\limits_{Q_w}|f_i(z)|^{\sigma_i}t^{\alpha_i}dz\bigg)^{\frac{q}{\sigma_i}}dw.$$

2. The measure $\mu$ satisfies the following Carleson-type condition:
$$\mu(\Delta_k)\leq C \eta_k^{m(n+1)\frac{p}{q}+\sum\limits_{i=1}^m\frac{p(n+1+\alpha_i)}{\sigma_i}},  k\geq 1.$$
\end{theorem}

\begin{theorem} (see \cite{3}). 
Let $0<p_i,\sigma_i<\infty, -1<\alpha_i<\infty,$ for $i=1,...,m.$ Let $\mu$ be a positive
Borel measure on $\mathbb{H}$. Then the following two conditions are equivalent:

1.  For any $m$-tuple $(f_1,..., f_m)$ of harmonic functions on $\mathbb{H}$ we have
$$\int\limits_{\mathbb{H}}\prod\limits_{i=1}^m|f_i(z)|^{p_i}d\mu(z)\leq C\prod\limits_{i=1}^m \int\limits_{\mathbb{H}}\bigg ( \int\limits_{Q_w}|f_i(z)|^{\sigma_i}dm_{\alpha_i}(z)\bigg)^{\frac{p_i}{\sigma_i}}dw.$$

2. The measure $\mu$ satisfies the following Carleson-type condition:
$$ \mu(\Delta_k)\leq C \eta_k^{m(n+1)+\sum\limits_{i=1}^m\frac{p_i(n+1+\alpha_i)}{\sigma_i}}, k\geq 1.$$
\end{theorem}

\begin{theorem} (see \cite{3}). 
Let $\mu$ be a positive Borel measure on an open proper subset $G$ of $\mathbb{R}^n$
 and let $(X_i,\|\cdot\|_{X_i})$ be a (quasi)-normed space of functions harmonic in $G, i=1,..,k.$ Let $d(x)=dist(x,\partial G), x\in G$ and $0<q_i<\infty, -1<\beta_i<\infty, i=1,..,k.$ Assume $$\sup\limits_{x\in G}|f_i(x)|^{q_i}d(x)^{\beta_i}\leq C\|f_i\|_{X_i}^{q_i}, f_i\in X_i, i=2,..,k,$$ and
 $$ \int\limits_{G}|f_1(x)|^{q_1}d\mu(x)\leq C\|f_1\|_{X_1}^{q_1}, f_1\in X_1.$$ Then
 $$ \int\limits_{G}\prod\limits_{i=1}^k|f_i(x)|^{q_i}d(x)^{\beta_2+...+\beta_k}d\mu(x)\leq C\prod\limits_{i=1}^k\|f_i\|_{X_i}^{q_i}, f_i\in X_i.$$
\end{theorem}

 \begin{theorem} (see \cite{3}). 
Let $0<p<q<\infty, \alpha>0$ and let $\mu$ be a positive Borel measure on $\mathbb{H}$.
Then the following two conditions are equivalent:

1. $$\int\limits_{\mathbb{H}}\bigg(\int\limits_{\mathbb{H}}\big(\frac{t^n}{|z-\overline{\omega}|^{2n}}\big)^{1+\alpha q}d\mu(z) \bigg)^{\frac{q}{q-p}}dm_{\alpha q n-1}(\omega)<\infty.$$

2.
$$\sum\limits_{k=1}^{\infty}\eta_k^{-(\alpha q n+n)\frac{p}{q-p}}\mu(\Delta_k)^{\frac{q}{q-p}}<\infty.$$
\end{theorem}

We define below new Herz type spaces of harmonic  functions in the unit ball similarly based on $r$-lattices we can define such type function spaces in upper half plane. Based on this definition we formulate several new sharp embedding theorems for such function spaces. Various related results in harmonic function spaces can be seen in papers of first author with M.~Arsenovic which we listed in references.

 Let $H_k$ be usual Whitney decomposition of upperhalf space $H$.

\begin{definition} Let $\mu$ be a positive Borel measure on $\mathbb{H}$ and let $0<p,q<\infty.$ We define
the space $K_q^p(\mu)$  as the space of $f\in h(\mathbb{H})$ such that
$$\|f\|^{q}_{K_q^p(\mu)}=\sum\limits_{-\infty}^{\infty}\bigg( \int\limits_{H_k}|f(z)|^p d\mu(z)\bigg)^{\frac{q}{p}}<\infty.$$
\end{definition}

 \begin{theorem} (see \cite{3}). 
Let $\alpha >-1, 0<p\leq q <\infty,$ and let $\mu$ be a positive Borel measure on $\mathbb{H}$.
Then the following two conditions are equivalent:

1. $A(p,\alpha)\hookrightarrow K_q^p(\mu).$

2. The measure $\mu$ satisfies the following Carleson type condition:
$$\frac{\mu(\Delta_j)}{|\Delta_j|^{1+\frac{\alpha}{n+1}}}\leq C, j\geq 1.$$
\end{theorem}

In the series of the following new sharp embedding theorems in higher dimension we collect namely several new embeddings in polydisk, unit ball, bounded strongly pseudoconvex domains obtained by the first author and his coauthors in recent years.
We refer for various notations and definitiins of basic objects related with these domains to references. To be more precise  we consider  below various analytic Herz type spaces, mixed norm spaces in product domains, Besov type function spaces in product domains.Proofs of some results below may be seen as some kind of modification of previously known sharp embedding  results obtained by other authors inmuch simpler stadard Bergman spaces in these domains. Some sharp embedding results below were proved by other authors in simpler domains.

Let $D=\{z: \rho (z)<0\}$ be a bounded strictly pseudoconvex domain in $C^n$
 with $C^{\infty}$  boundary. We assume that the strictly plurisubharmonic function $\rho$ is of class $C^{\infty}$
in a neighbourhood of $\overline{D},$ that $-1\leq \rho(z)<0, z\in D, |\partial\rho|\geq c_0>0$ for $|\rho|\leq r_0.$

Denote by $H(D)$ the space of all analytic functions on $D$. Let also $A^{p,q}_{\delta,k}=\{f\in H(D): \|f\|_{p,q,\delta,k}<\infty\},$ where

$$\|f\|_{p,q,\delta,k}=\bigg( \sum\limits_{|\alpha|\leq k}\int\limits_0^{r_0}\bigg( \int\limits_{\partial D_r}|D^{\alpha}f|^p d\sigma_r\bigg)^{\frac{q}{p}}r^{\delta\frac{q}{p}-1}dr \bigg)^{\frac{1}{q}}
$$
be the mixed norm space in $D.$ Here $D_r=\{z\in \mathbb{C}^n : \rho(z)<(-r)\}, \partial D_r$ is boundary $d\sigma_r$ the normalized surface measure on $\partial D_r$ and by $dr$ we denote the normalized volume element on
interval from 0 to $r, 0<p<\infty, 0<q\leq \infty, \delta>0, k=0,1,2...,$ and
$$\|f\|_{p,\infty,\delta,q}=\sup  \bigg \{ ( \sum\limits_{|\alpha|\leq k}(r^\delta\int\limits_{\partial D_r}|D^{\alpha}f|^p)d\sigma_r \bigg)^{\frac{1}{p}} : 0<r<r_0 \bigg \} $$
(for $p, q < 1$ it is quazinorm), where $D^{\alpha}$ is a derivative of $f.$\\
For $p=q$ we  have
$$\|f\|_{p,\delta,k}=\bigg(\sum\limits_{|\alpha|\leq k}\int\limits_{D}|D^{\alpha}f|^p(-\rho)^{\delta-1}d\zeta\bigg)^{\frac{1}{p}}; \delta>0, k\geq 0.$$

\begin{theorem} (see \cite{4}). 
Let $D$ be a bounded strictly pseudoconvex domains with smooth boundary. Let $\mu$
be positive finite Borel measure on $D, f\in H(D).$ Let $\{a_k \}$ be r-lattice. Assume $q<p$ or $q=p, r\leq p.$ Then we have
$$\bigg(\int\limits_D |f(z)|^p d\mu(z) \bigg)^{\frac{1}{p}}\leq c_0\|f\|_{A_{v,0}^{q,\tau}} $$
if and only if $\mu(B_D(a,r))\leq c_1 \delta^{\frac{np}{q}+\frac{vp}{q}}(a), a\in D,$ or if and only if $\mu(B_D(a_k,r))\leq c_2 \delta^{\frac{np}{q}+\frac{vp}{q}}(a_k)$ for $k=1,2,...,$ and for some constants $c_1,c_2$ for all $v< \big(\frac{2q}{p}-1 \big)n+2\frac{q}{p},$ if $q<p$ and for all $v<(n+2)\frac{p}{r},$ if $q=p,r\leq p.$
\end{theorem}

Let
$$U^n=\{ z=(z_1,...,z_n): |z_j|<1, 1\leq j\leq n\}$$
be the unit polydisk of $n$-dimensional complex space $\mathbb{C}^n, T^n$ be the Shilov boundary of $U^n, \vec{p}=(p_1,...,p_n), 0<p_j<\infty,$ $j=1,...,n, \vec{w}(t)=(w_1(t),...,w_n(t)), t\in (0,1),$ where $w_j(t)$ are positive integrable functions on $(0,1).$ We denote by $A^{\vec{p}}(\vec{w})$ the set
of all holomorphic functions in $U^n$ for which
$$\|f\|_{A^{\vec{p}}(\vec{w})}=\bigg ( \int\limits_U \bigg [...\big ( \int\limits_U |f(z_1,...,z_n)|^{p_1}\omega_1(1-|z_1|)dm_2(z_1)\big)^\frac{p_2}{p_1} ...\bigg  ]^{\frac{p_n}{p_n-1}}\omega_n(1-|z_n|)dm_2(z_n)\bigg )^{\frac{1}{p_n}}<\infty,
$$
where $m_2$ is planar Lebesgue measure on $U:=U^1.$ Assume further $\vec{\mu}=(\mu_1,...,\mu_n),$ where $\mu_j$
is the Borel nonnegative finite measure on $U, L^{\vec{p}}(\vec{\mu})$ is related
space with mixed norm that is, the space of all measurable functions $U^n$
for which
$$\|f\|_{L^{\vec{p}}(\vec{\mu})}=\bigg( \int\limits_U \bigg[...\big( \int\limits_U |f(\zeta_1,...,\zeta_n)|^{p_1}d\mu_2(\zeta_1)\big)^{\frac{p_2}{p_1}} ...\bigg]^{\frac{p_n}{p_n-1}}d\mu_n(\zeta_n) \bigg)^{\frac{1}{p_n}}<\infty,
$$
$0<p_i\leq\infty, i=1,...,n,$ with usual modification for $p_i=\infty.$

Let $S$ be the set of all measurable and positive functions of $L^1(0,1)$ or which there
exists numbers $M_\omega, m_\omega, q_\omega$ with $m_\omega, q_\omega\in(0,1]$ that is
$$m_\omega\leq \frac{\omega(\lambda r)}{\omega(r)}\leq M_\omega,  r\in (0,1], \lambda\in [q_\omega].$$

Let $L^{\vec{p}}(\vec{w},\vec{\nu})=$
$$=\bigg \{ f \in L_{loc}^1(B_k^n): \bigg ( \int\limits_{B}...\big( \int \limits_{B}|f(z_1,...,z_n)|^{p_1}(w_1(1-|z_1|)d\nu_1(z_1))^{\frac{p_2}{p_1}} \big )... (w_n(1-|z_n|)d\nu_n(z_n))^{\frac{1}{p_n}}\bigg) < \infty \bigg \}, $$
$0<p_i<\infty,$ $ i=1,...,n,$ $ \nu_j, j=1,...,n,$ be the normalized Lebegues measures on $B_k, A^{\vec{p}}(\overrightarrow{w})=L^{\vec{p}}(\vec{w})\cap H(B_k^n).$ Replacing $w_jd\nu_j$ by $d\mu_j$ we define similarly the new
general space $L^{p_1,...,p_n}(\mu_1,...,\mu_n).$

Let $$\tilde{\mathcal{D}}^{\vec{\alpha}}f(\vec{z})=\sum\limits_{|\vec{k}|\leq 0}\prod_{j=1}^m(|k_j|+1)^{\alpha_j}a_{k_1,...,k_m}z_1^{k_1}...z_m^{k_m},$$
where $\sum_{|\vec{k}|\geq 0}$ means $\sum_{k_1\geq 0}...\sum_{k_m\geq 0}.$

We extend in a natural way (as in polydisk case) the definition of differential
operator $\widetilde{D}^m$  to differential operators acting on analytic functions defined on product
domains for all real $\alpha_j, j=1,...,m.$

In the following two theorems we show new interesting sharp embeddings in so called new mixed norm spaces of analytic functions in product domains which were introduced for the first time  F. Shamoyan in his interesting paper \cite{d7}.
 
Actually, as far as we know, both theorems
below belong to him but he showed them in context of less general polydisk.We proved them later in general polyball case. In case of function spaces in $R^n$ such mixed norm interesting function spaces were well studied by S. Nikolski, I. Ilin (see, for example, \cite{k8}) and many other, we refer to \cite{d7} for these books and papers concerning the case of function spaces in $R^n$ and for history of problem, concerning these interesting spaces.
These function spaces of analytic functions are
completely new and we note that
many problems are open in this research area.
We refer the reader to pay attention to them and leave them to interested readers.

\begin{theorem} (see \cite{5}). 
Let $\{a_k\}$ be $r-$lattice of $B_k$. Let $\vec{\omega}=(\omega_1,...,\omega_n), \vec{\mu}=(\mu_1,...,\mu_n),\omega_j\in S, m=(m_1,...,m_n)\in \mathbb{Z}^n_{+}, \vec{p}=(p_1,...,p_n),\vec{q}=(q_1,...,q_n)\in \mathbb{R}^n_{+}$ with $0<p_j\leq q_j, j=1,...,n.$ Then the following assertions are equivalent:

1. $\|\tilde{D}^m(f)\|_{L^{\vec{q}}(\vec{\mu})}\leq C(\vec{\mu})\|f\|_{A^{\vec{p}}(\vec{\omega})};$

2. $\mu_j(D(a_k,r))\leq c(1-|a_k|)^{(n+1)\frac{q_j}{p_j}+m_jq_j}[\omega_j(1-|a_k|)]^{\frac{q_j}{p_j}}, j=1,..,n; k=0,1,2,...$
\end{theorem}

In the case of measures $\nu$ defined on $B_k^n=B_k\times...\times B_k$ there is following result.

\begin{theorem} (see \cite{5}).  
Let $p_j<q<+\infty, \tilde{\nu}$ be the Borel nonnegative measure on $B_k^n, $ $\vec{\omega}=(\omega_1,...,\omega_n), \omega_j\in S, j=1,...,n, m=(m_1,...,m_n)\in\mathbb{Z}^{n}_{+}.$ Then the following assertion are equivalent:

1. $\big( \int\limits_{B^n_k}|\tilde{D}^mf(z)|^{q}d\tilde{\nu}(z)\big)^{\frac{1}{q}}\leq C\|f\|_{A_{\vec{p}}(\vec{\omega})};$

2. $\tilde{\nu}(D(a_{k_1},r)\times ...\times D(a_{k_n},r))\leq c\prod\limits_{j=1}^n (1-|a_{k_{j}}|)^{(n+1)\frac{q}{p_j}+m_j q}[\omega_j(1-|a_{k_j}|)]^{\frac{q}{p_j}}.$
\end{theorem}

Let $A$ be a Banach space
of analytic functions on a domain $D\subset\mathbb{C}^n.$  Given $p > 1$, a finite positive Borel
measure $\mu$ on $D$ is a Carleson measure of $A$ (for $p$) if there is a continuous inclusion $A\hookrightarrow A^{p}(\mu),$ that is, there exists a constant $C > 0$ such that
$$\int\limits_{D}|f|^p d \mu\leq C\|f\|^{p}_A, \forall f\in A.
$$
A finite positive Borel measure $\mu$ is a Carleson measure of Hardy space $H^p(\mathbb{D}),$ where $ \mathbb{D}$ is the unit disk in $\mathbb{C}$ if and only if there exists a constant $C > 0$ such that $\mu(S_{\theta_0,h})\leq Ch$ for all sets
$$S_{\theta_0,h}=\{re^{i\theta}\in \mathbb{D}:1-h\leq r<1, |\theta-\theta_0|<h\} .$$
The set of Carleson measures of $H^p(\mathbb{D})$  does not depend on p.

Let $z_0\in D$ and $0<r<1, $ let $B_D(z_0,r)$ denote the ball of center $z_0$ and radius $\frac{2}{2}(\log(1+r)-\log(1-r))$ for the Kobayashi
distance $k_D$ of $D.$ It is known for $D$ strongly pseudoconvex  that a finite
positive measure $\mu$ is a Carleson measure of $A^p(D)$ for $p$ if and only if for some (and hence all) $0<r<1$ there is a constant $C_r>0$ such that
$$\mu(B_D(z_0,r))\leq C_r \nu(B_D(z_0,r)), \forall z_0\in D,
$$
where $\nu$ is a standard Lebegues measure on a domain.

We say that a finite positive Borel measure $\mu$ is a (geometric) $\theta$-Carleson measure
if for some (and hence all) $0<r<1$ there is a constant $C_r>0$ such that
$$\mu(B_D(z_0,r))\leq c_r \nu(B_D(z_0,r))^{\theta}, \forall z_0\in D.$$

We define new analytic  Herz type function spaces of several complex variables in general $D$ domain below where so called $r$-lattices exists. It is well known that the existence of these $r$-lattices in various domains in $C^n$ and it is various nice properties is a base for many new results in particular embeddinds in analytic spaces in such domains in complex function spaces of several variables.
For such Herz type function spaces of several complex variables we below provide a series of new sharp embedding theorems provided in recent year by first author and his coauthors.
We can also define similar analytic Herz spaces replacing summation by integration below.
These another type of new Herz spaces is a new research area and similar theorems (embeddings) can be provided probably there.
We leave these problems to interested readers.

In the series of the following embedding theorems we provide sharp embeddings in context of $D$ bounded strongly pseudoconvex domains.

\begin{definition}
Let $\mu$ be a positive Borel measure in $D, 0<p,q<\infty, s>-1.$ Fix $r\in(0, \infty),$ and an $r$-lattice $\{a_k\}_{k=1}^{\infty}.$ The analytic space $A(p,q,d\mu)$ is the space of all holomorphic functions $f$ such that
$$\|f\|^{q}_{A(p,q,d\mu)}=\sum\limits_{k=1}^{\infty}\bigg(\int\limits_{B(a_k,r)}|f(z)|^p d\mu(z)\bigg)^{\frac{q}{p}}<\infty.
$$
\end{definition}

If $ d\mu=\delta^s(z)d\nu(z),$ then we will denote by $A(p,q,s)$ the space $A(p,q,d\mu).$

\begin{theorem} (see \cite{6}). 
Let $\mu$ be a positive Borel measure on $D$ and $\{a_k\}$ be a Bergman
sampling sequence forming an $r$-lattice. Let $\alpha>-1, f_i\in H(D), 0<p_i,q_i<\infty, i=1,...,m$ so that $\sum\limits_{i=1}^m q_i^{-1}=1.$ Then $$\int\limits_D \prod\limits_{i=1}^m|f_i(z)|^{p_i}d \mu(z)\leq C\prod \limits_{i=1}^m \bigg [ \sum\limits_{k=1}^{\infty} \big ( \int\limits_{B(a_k,r)}|f_i(z)|^{p_i}\delta^{\alpha}(z)d\nu(z)\big)^{q_i}\bigg]^{\frac{1}{q_i}}
$$
if and only if
$$\mu(B_D(a_k,r))\leq C\delta^{m(\alpha+2n/r)}(a_k)
$$
for every $k=1,2,3,..., r>0.$
\end{theorem}

We denote by $d\nu_\alpha$ standard weighted Lebegues measure in theorems below in various domains (see \cite{6}).

\begin{theorem} (see \cite{6}). 
Let $0<q,s<\infty, q\geq s, \alpha>-1.$ Let $\{a_k\}_{k=1}^{\infty}$ be a sequence
forming an $r$-lattice in $D$. Let $\mu$ be a positive Borel measure in $D$. Then
$$\int\limits_{D}|f(z)|^q d\mu(z)\leq C\int\limits_D \bigg(\int\limits_{B_D(z,r)}|f(w)|^s d\nu_\alpha(w) \bigg)^{\frac{q}{s}}d\nu(z)
$$
if and only if
$$\mu(B_D(a_k,r))\leq C(\delta(a_k))^{q((\alpha+2n/r)/s+2n/(rq))}
$$
for some constant $C, C>0, k\in \mathbb{N}.$
\end{theorem}

\begin{theorem} (see \cite{6}). 
Let $0<q,p<\infty, 0<s\leq p<\infty, \beta>-1, \tilde{\beta}=\beta+nr^{-1}.$  Let $\mu$ be a positive Borel measure on $D$. Then we have the assertion
$$\|f\|_{A(p,q,d\mu)}\leq C\|f\|_{A^s_{\tilde{\beta}}}
$$
if and only if
$$\mu(B_D(a_k,r))\leq C(\delta(a_k))^{\frac {q(2n/r+\beta)}{s}}.
$$
\end{theorem}

Let $A$ be a Banach space of analytic functions on a domain $D\in \mathbb{C}^n$. Given $p \geq 1$, a finite positive Borel
measure $\mu$ on $D$ is a Carleson measure of $A$ (for $p$) if there is a continuous inclusion $A\hookrightarrow L^p(\mu)$, that is, there exists a constant $C > 0$ such that
$$\int\limits_D |f|^p d\mu \leq C\|f\|^p_A, \forall f\in A.$$
A finite positive Borel measure $\mu$ is a Carleson measure of Hardy space $H^p(\mathbb{D}),$
where $\mathbb{D}$ is the unit disk in $\mathbb{C}$  if and only if there exists a constant $C > 0$ such that $\mu(S_{\theta_0,h})\leq Ch$ for all sets
$$ S_{\theta_0,h}=\{r e^{i\varphi} \in \mathbb{D}: |1-h|\leq r<1, |\theta-\theta_0|<h \}.$$
The set of Carleson measures of $H^p(\mathbb{D})$  does not depend on $p$.

We define Herz spaced and formulate below complete analogues of theorems above but in context of bounded strongly pseudoconvex domains.

\begin{definition} Let $\mu$ be a positive Borel measure in $D, 0<p, q<\infty, s>-1.$ Fix $r\in(0,\infty)$ and $r$-lattice $\{a_k\}_{k=1}^{\infty}.$ The analytic space $A(p,q,d\mu)$ is the space
of all holomorphic functions $f$ such that
$$ \|f\|^{q}_A(p,q,d\mu)=\sum\limits_{k=1}^{\infty} \bigg(\int\limits_{B(a_k,r)}|f(z)|^p d\mu(z) \bigg)^{\frac{q}{p}}< \infty.$$
\end{definition}

If $d\mu=\delta^s(z)d\nu(z),$ then we will denote by $A(p,q,s)$ the space $A(p,q,d\mu).$ This is
a Banach space for $\min(p,q) \geq 1.$ It is clear that $A(p,p,s)=A^p_{\bar{s}}, \bar{s}=s+nr^{-1}.$

Let $\mu$ be a positive Borel measure in tube $D$ and let $\{a_k\}_{k\in \mathbb{N}}$ be a sequence such that $B_D(a_k,r)$ is an $r$-lattice for tubular domain $D$ in $\mathbb{C}^n.$ Let $X$  be
a quasinormed subspace of $H(D)$ and $p,q\in (0, \infty)$. Describe all positive Borel measures such that
$$ \|f\|_{A(p,q,d\mu)}\leq C \|f\|_X.$$

We provide below complete analogues of previous three theorems in tube in context of bounded strongly pseudoconvex domains. Formulations of these theorems are very similar.

\begin{theorem} (see \cite{7}). 
Let $\mu$ be a positive Borel measure on $D$ and $\{a_k\}$ be a Bergman
sampling sequence forming an $r$-lattice. Let $\alpha >-1, f_i\in H(D), 0<p_i, q_i<\infty, i=1,...m,$ so that $\sum\limits_{i=1}^m q_i^{-1}=1.$ Then
$$\int\limits_D \prod\limits_{i=1}^m |f_i(z)|^{p_i} d\mu(z)\leq C \prod \limits_{i=1}^m \bigg [ \sum \limits_{k=1}^{\infty} \bigg( \int\limits_{B(a_k,r)} |f_i(z)|^p_i \delta^{\alpha}(z)d\nu(z) \bigg)^{q_i}\ \bigg]^{\frac{1}{q_i}}
$$
if and only if
$$\mu(B_D(a_k,r))\leq C \delta^{m(a+\frac{2n}{r})}(a_k)$$
for every $k=1,2,..., r>0.$
\end{theorem}

\begin{theorem} (see \cite{7}). 
Let $0<q, s<\infty, q\geq s, \alpha>-1.$ Let $\{a_k\}_{k=1}^{\infty}$-be a sequence
forming an $r$-lattice in $D$. Let $\mu$ be a positive Borel measure in $D$. Then
$$\int\limits_D |f(z)|^q d\mu(z)\leq C \int\limits_D \bigg( \int\limits_{B_D(z,r)} |f(w)|^s d\nu_{\alpha}(w) \bigg )^{\frac{q}{s}}d\nu(z)
$$
if and only if
$$ \mu(B_D(a_k,r))\leq C (\delta(a_k))^{q((\alpha+2n/n)/s+2n/(rq))}$$
for some constant $C, C > 0, k\in \mathbb{N}$.
\end{theorem}

\begin{theorem} (see \cite{7}). 
Let $0<q, p<\infty, 0<s\leq p< \infty, \beta>-1, \tilde{\beta}=\beta+nr^{-1}.$ Let $\mu$ be a positive Borel measure on $D$. Then we have the assertion
$$\|f\|_{A(p,q,d\mu)}\leq C \|f\|_{A_{\tilde{\beta}}^s}$$
if and only if
$$\mu(B_D(a_k,r))\leq C(\delta(a_k))^{q(2n/r+\beta)/s}.$$
\end{theorem}

Let $D=\{z: \rho(z)<0\}$ be a bounded strictly pseudoconvex domain of $\mathbb{C}^n$ with $\mathbb{C}^{\infty}$ boundary.We assume that the strictly plurisubharmonic function $\rho$ is of class $\mathbb{C}^{\infty}$ in a neighbourhood of $\bar{D},$ that, $-1\leq \rho(z)<0, z\in D, |\partial\rho|\geq c_0>0$ for $|\rho|<r_0.$

Denote by $\mathcal{O}(D)$ or $H(D)$ the spaces of all analytic functions on $D.$ Let also
$A^{p,q}_{\delta,k}=\{f\in H(D): \|f\|_{p,q,\delta,k}<\infty\},
$ where
$$ \|f\|_{p,q,\delta,k}=\bigg ( \sum\limits_{|\alpha|\leq k} \int\limits_{0}^{\rho_0} \big(\int\limits_{\partial D_r} |D^{\alpha}f|^p d\sigma_r \big)^{\frac{q}{p}} r^{\delta \frac{q}{r}-1} dr\bigg)^{\frac{1}{q}}$$
is the mixed norm on $D.$ Here $D_r=\{z\in \mathbb{C}^n: \rho(z)<(-r)\}, \partial D_r$ is
its boundary, $d\sigma_r$ is the normalized surface measure on $\partial D_r,$ and by $d r$ we
denote the normalized volume element on $D, 0<p<\infty, 0<q\leq \infty, \delta>0, $ $ k=0,1,2,...,$  and
$$\|f\|_{p,\infty,\delta,k}=\sup \bigg \{ \bigg(\sum\limits_{|\alpha|\leq k} (r^\delta) \int\limits_{\partial D_r}|D^{\alpha}f|^p d\sigma_r \bigg)^{\frac{1}{p}} :0<r<r_0 \bigg \},
$$
where $D^{\alpha}$ is a derivative of $f.$

For $p=q$ we have $$ \|f\|_{p,\delta,k}=\bigg ( \sum\limits_{|\alpha|\leq k} \int\limits_{D}|D^{\alpha}f|^p (-\rho)^{\delta}d \zeta\bigg)^\frac{1}{p}; \delta > -1, k\geq 0.$$
For $p = q, k = 0$ we get the usual Bergman spaces $A_\delta^p(D).$

We use for convenience the following notations below: $Q_w=B(w,r), w\in D, \{a_k\}$ is a fixed $r$-lattice in $D,$ and $\Delta_k^{*}=|B^*(a_k,R)|, R=\frac{1+r}{2}, \Delta_k=|B(a_k,r)|.$
We refer for $B(z,r)$ to \cite {8}. 

In our theorems 25 and 28 below we assume that
$$ \int_{B(z,r)}|K_{\tau}(\tilde{z},w)|^p \delta^{\alpha}(\tilde{z})dv(\tilde{z})\leq C|K_{\tilde{\tau}}(w,z)|,$$
$\tilde{\tau}=\tau p-\alpha -n-1, 0<p<\infty, \alpha>-1, \tau>0, w,z \in D, \tau$ is large enough.

In the following embedding theorems and estimates we provide new sharp embedding theorems and estimates for analytic function spaces
Related with the so called new trace function.
Such type sharp results and estimates are new and they may be valid in many other domains such as bounded symmetric domains with very similar proof or in Siegel domains of second type.
Proofs in all domains namely in tubular domains, bounded pseudoconvex domains and in harmonic function spaces are similar they were obtained by the first author and his various coauthors in recent years.

We define Trace operator in any $D$ domain in $C^n,$ as a map $f(z_1,...z_m) \rightarrow f (z,...,z), z, z_j \in D.$

\begin{theorem} (see \cite{8}). 
Let $0<p<\infty$ and let, for $i=1,...,m, 0<q_i, \sigma_i<\infty, \alpha_i>-1.$  Assume that $q_i \leq p$ for $i=1,...,m.$ Let $\mu$ be a positive Borel measure on $D$. Then the following two conditions are equivalent:

1. For any analytic function $f(z_1,...,z_m)$ on $D^m$ that splits into a product
of analytic functions $f_i(z_i)\in H(D),$  i.e.  $f(z_1,...,z_m)=\prod_{j=1}^m f_j(z_j)$ we have
$$\bigg(\int_D |Trf(z)|^p d\mu(z)\bigg)^{\frac{1}{p}}\leq C \prod\limits_{i=1}^m \bigg( \int_D \big(  \int_{Q_w} |f_i(z)|^{\sigma_i} dv_{\alpha_i}(z)\big)^{\frac{q_i}{\sigma_i}} dv(w)\bigg)^{\frac{1}{q_i}}.$$

2. The measure $\mu$ satisfies the following Carleson-type condition:
$$\mu(\Delta_k)\leq C\delta(a_k)^{p\sum\limits_{i=1}^m (\frac{n+1+\alpha_i}{\sigma_i}+\frac{n+1}{q_i})}, k\geq 1.$$
\end{theorem}

\begin{theorem} (see \cite{8}). 
Let $0<p<\infty, s_j>-1,j=1,...,m.$ Let $\mu$ be a positive Borel measure on $D.$ Then the following are equivalent:

1. $\mu(\Delta_k)\leq C|\Delta_k|^{m+\frac{1}{n}\sum\limits_{j=1}^m s_j}, k\geq 1;$

2.  $(Trace)(A^{p}_{\vec{s}})\hookrightarrow L^p(D,d\mu);$

3.$\int_D \prod\limits_{j=1}^m |f_j(z)|^p d\mu(z)\leq c\prod\limits_{j=1}^m \|f_j\|^p_{A^p_{s_j}}, f_j\in A^p_{s_j}, 1\leq j\leq m.$
\end{theorem}

\begin{theorem} (see \cite{8}). 
Let $\alpha>-1, f_i\in H(D^t),$ for $i=1,...,m.$ Let $0<p_i,q_i<\infty, \sum\limits_{i=1}^m \frac{p_i}{q_i}=1$, $\beta_i=\frac{(n+1+\alpha)q_i}{tmp_i}-(n+1)>-1, i\in [1,m]. $ Then

$$ \int_D |Trf_1(w)|^{p_1}...|Trf_m(w)|^p_m \delta(w)^{\alpha}dv(w)\leq$$
$$\leq C\prod\limits_{i=1}^m \bigg(\int_D...\int_D |f_i(w_1,...,w_t)|^{q_i} \bigg( \prod\limits_{j=1}^t \delta(w_j)^{\beta_j}\bigg) dv(w_1)...dv(w_t) \bigg)^{\frac{p_i}{q_i}}.$$
\end{theorem}

We put an additional assumption for the following sharp theorem (Carleson type theorem). Let $f\in H(D), \tau \geq 0, $ $ \alpha>-1,$  and assume that
$$|\delta(a_k)|^{n+1} \bigg(\int_{B^*(a_k,r)}|f(z)|dv_\alpha(z)\bigg)^{\tau} \leq c \int_{B^*(a_k,r)} \bigg( \int_{B(w,r)}|f(z)|dv_{\alpha}(z)\bigg)^{\tau}dv(w).$$

\begin{theorem} (see \cite{8}). 
Let $0<p_i,\sigma_i<\infty, -1<\alpha_i \leq \infty$ for $i=1,...,m.$  Let $\mu$ be a positive Borel measure on $D$. Then the following two conditions are equivalent:

1. For any $m$-tuple $(f_1,...,f_m)$  of analytic functions on $D$ we have
$$ \int_D \prod\limits_{i=1}^m |f_i(z)|^{p_i}d \mu(z) \leq C\prod\limits_{i=1}^m \int_D \bigg( \int_{Q_w} |f_i(z)|^{\sigma_i}dv_{\alpha_i}(z)\bigg)^{\frac{p_i}{\sigma_i}}dv(w).$$

2. The measure $\mu$ satisfies a Carleson type condition:
 $$ \mu(\Delta_k)\leq C\delta(a_k)^{m(n+1)+\sum\limits_{i=1}^{m}\frac{p_i(n+1+\alpha_i)}{\sigma_i}}, k\geq 1. $$
\end{theorem}

 Let $T_\Omega$ is a tube domain over a cone, $z_0\in D$ and $r\in (0,1)$. We denote
by $B(z_0,r)$ the Bergman ball of center $z_0$ and radius $r$. We denote by $dV$ the
normalized Lebegues measure on $T_\Omega$.
Let $T_\Omega=V=i\Omega$ be the tube domain over an  irreducible symmetric cone $\Omega$
in the complexification $V^{\mathbb{C}}$ of an $n$-dimensional Euclidean space $V$. We denote
the rank of the cone $\Omega$ by $r$ and by $\Delta$ the determinant function on $V$.

$H(T_\Omega)$ denotes the space of all holomorphic functions on $T_\Omega$. Also, we denote
$m$-products of tubes by $T_\Omega^m=T_\Omega\times ...\times T_\Omega.$ The space of all analytic function
on this new product domain which are analytic by each variable separately will
be denoted by $H(T_\Omega^m).$

For $\tau \in \mathbb{R}_+$ and the associated determinant function $\Delta (x)$ we set
$$A_\tau^{\infty}(T_\Omega)=\bigg \{ F\in H(T_\Omega): \|F\|_{A_\tau^{\infty}}= \sup \limits_{x+iy\in T_\Omega} |F(x+iy)|\Delta^{\tau}(y)<\infty\bigg\}.
$$
For $1\leq p,q< +\infty, \nu\in \mathbb{R}$ and $\nu> \frac{n}{r}-1$ we denote by $A_{\nu}^{p,q}(T_\Omega)$ the
mixed-norm weighted Bergman space consisting of analytic functions $f$ in $T_\Omega$ that
$$\|F\|_{A_{\nu}^{p,q}}= \bigg (\int\limits_\Omega \bigg( \int\limits_ V |F(x+iy)|^p dx \bigg)^{\frac{q}{p}} \Delta^{\nu}(y) \frac{dy}{ \Delta(y)^{\frac{n}{r}}}  \bigg)^{\frac{1}{q}}<\infty.
$$
This is also a Banach space. It is known the $A_\nu^{p,q}(T_\Omega) $ space is nontrivial if
and only if $\nu >\frac{n}{r}-1$ and we will assume this. When $p=q$ we write
$$A_\nu^{p,q}(T_\Omega)=A_\nu^p (T_\Omega). $$

In theorems 29 and 31 we consider new type of sharp  embeddings in analytic spaces on product domains related with the so called trace operator.in tubular domains over symmetric cones. Complete analogues of these results in harmonic function spaces of several variables were obtained by the first author with M. Arsenovic. These results are also valid in bounded strongly pseudoconvex domains with the smooth boundary. We note
similar type results may be valid also in various other domains.

Some sharp embedding results related with sharp embedding  in analytic spaces related with the trace operator which we see below can be seen also in a paper of the first author with E. Tomashevskaya (see \cite {d13}). Our next embedding results we formulate in context of tube domains though these resultscan be valid with the same or similar proofs in context of other domains which we discussed shortly above also.

\begin{theorem} (see \cite{8}). 
Let $0<p<\infty$ and let, for $i=1,...,m, 0<q_i,\sigma_i<\infty, \alpha_i>-1.$ Assume that $q_i\leq p$ for $i=1,...,m.$  Let $\mu$ be a positive Borel measure on $T_\Omega$. Then the following two conditions are equivalent:

1. For any analytic functions $f(z_1,...,z_m)$ on $T_\Omega^m$ that splits into a product
of analytic functions $f_i(z_i)\in H(T_\Omega), i=1,...,m,$ i.e. $f(z_1,...,z_m)=\prod\limits_{j=1}^m f_j(z_j)$ we have
$$\bigg( \int\limits_{T_\Omega}|Tr f(z) |^p d\mu(z)\bigg)^{\frac{1}{p}}\leq C\prod\limits_{i=1}^m \bigg( \int\limits_{T_\Omega}\big(  \int\limits_{\tilde{B}(w,r)} |f_i(z)|^{\sigma_i} \Delta^{\alpha_i} Im (z)dV(z)\big )^{\frac{\sigma_i}{q_i}} dV(w)\bigg)^{\frac{1}{q_i}}.
$$

2. The measure $\mu$ satisfies the following Carleson-type condition:
$$\mu(\tilde{B}(a_k,r))\leq C\Delta(Im (a_k))^{p\sum\limits_{i=1}^m\big ( \frac{\frac{2n}{r}+\alpha_i}{\sigma_i}+\frac{\frac{2n}{r}}{q_i} \big)}, k\geq 1.
$$
\end{theorem}

We assume further that
$$|\tilde{B}(a_k,r)|^{\frac{2n}{r}} \bigg( \int\limits_{\tilde{B}^*(a_k,r)}|f(z)|dV_{\alpha}(z) \bigg)^{\tau}\leq C\int\limits_{\tilde{B}^*(a_k,r)} \bigg(\int\limits_{\tilde{B}(w,r)}|f(z)|dV_{\alpha}(z)\bigg)^{\tau}dV(w),
$$
for any $\{a_k\}$ $r$-lattice in $T_\Omega$, for any Bergman ball $\tilde{B}(w,r)$ and any $f, f\in H(T_\Omega)$ and for $dV_{\alpha}(z)=\Delta^{\alpha}(Im z), $ $\alpha>-1, $ $ 0<\tau<\infty , $ $ z\in T_\Omega,$ where $\tilde{B}^*(a_k,r)$ is an enlarged dyadic Bergman ball.

\begin{theorem} (see \cite{8}). 
Let $0<p_i,\sigma_i<\infty, -1<\alpha_i\leq \infty$ for $i=1,...,m.$ Let
$\mu$ be a positive Borel measure on $T_\Omega$. Then the following two conditions are
equivalent:

1. For any $m$-tuple $(f_1,...,f_m)$ of analytic functions on $T_\Omega$ we have
$$\int\limits_{T_\Omega} \prod\limits_{i=1}^m |f_i(z)|^{p_i} d\mu(z)\leq C \prod\limits_{i=1}^m \int\limits_{T_\Omega} \bigg( \int\limits_{\tilde{B}(w,r)}|f_i(z)|^{\sigma_i} \Delta^{\alpha_i} Im(z) dV(z)\bigg)^{\frac{p_i}{\sigma_i}}dV(w).
$$

2. The measure $\mu$ satisfies a Carleson type condition:
$$\mu(\tilde{B}(a_k,r))\leq C\Delta(Im(a_k))^{m\frac{2n}{r}+\sum\limits_{i=1}^m\frac{p_i(\frac{2n}{r}+\alpha_i)}{\sigma_i}}, k\geq 1.$$
\end{theorem}

Let further for $m>1, H(T_\Omega^m)$ be the space of all analytic functions on $T_\Omega^m=T_\Omega\times...\times T_\Omega, m\in \mathbb{N}.$ Let
$$A_{\vec{s}}^p=\{f\in H(T_\Omega^m): \int\limits_{T_\Omega} ...\int\limits_{T_\Omega} |f(z_1,...,z_m)|^p \prod\limits_{j=1}^m \Delta^{s_j-\frac{n}{r}}(Im z)dV(z)<\infty\},
$$
where $s=(s_1,...,s_m), 0<p<\infty, s_j>\frac{n}{r}-1, j=1,...,m.$ Let also further $L^p(T_\Omega, d\mu)$ be a space of all measurable functions on $T_\Omega$ so that
$$\int\limits_{T_\Omega}|f(z)|^p d\mu(z)<\infty,$$
for any positive Borel measure $\mu$ on $T_\Omega$.

We wish to remark that a version of the following sharp embedding theorem which we formulate below in context of the tubular domains over symmetruc cone may be also valid in various other domains with the same or similar  proofs. We omit all details here. Leaving them to our interested readers.

\begin{theorem} (see \cite{8}). 
Let $0<p<\infty, s_j>-1, j=1,...,m.$ Let $\mu$ be a positive Borel measure on $T_\Omega$. Then the following are equivalent:

1. $\mu(\Delta_k) \leq C|\Delta_k|^{m+(\frac{2n}{r})^{-1}\sum_{j=1}^m s_j}, k\geq 1;$

2.$(Trace)(A_{\vec{s}}^p)\hookrightarrow L_p(T_\Omega, d\mu);$

3.$\int\limits_{T_\Omega}\prod\limits_{j=1}^m |f_j(z)|^p d\mu(z) \leq c\prod\limits_{j=1}^m \|f_j\|^p_{A_{s_j}^p}, f_j\in A_{s_j}^p, 1\leq j\leq m.$
\end{theorem}

\section{Final discussions and comments}

We refer to paper J. Ortega and J. Fabrega \cite{d12} for various interesting new embeddings in Bergman type and mixed norm analytic function spaces in bounded strongly pseudoconvex domains.
We pose as a problem to define similar type more complicated multifunctional Bergman type and mixed norm analytic function spaces and to  extend these embedding results to  multifunctional analytic Bergman type  function spaces in bounded strongly pseudoconvex domains with the smooth boundary by modification of proofs from mentioned paper.

Many rather interesting open  problems related with embeddings can be posed  also in  analytic area Nevanlinna spaces of several complex variables. We leave this task to interested readers. Note first various types of such type Nevanlinna type analytic function spaces must be defined in several complex variables and then direct relations (embeddings) between such type function spaces must be studied. Even in one functional analytic spaces of several complex variables various open problems in this research area can be posed.

Sharp new interesting embeddings for analytic area Nevanlinna type spaces of one complex variables were obtained recently in \cite{d332}. We refer the reader to this paper for the story of problem and other interesting results in this direction.
The natural interesting new problem here to define similar type spaces of area Nevanlinna type in several complex variables in various complicated domains and to find similar type sharp embedding results in several complex variables for analytic function spaces of area Nevanlinna type. We pose this as a problem for interested readers.

It will be very interesting to develop the techniques used by Igor Verbitsky which we shortly mentioned above  in context of the unit disk to spaces with the mixed norm of analytic functions of similar type but in context of the unit polydisk and even more difficult domains. For this case in particular analytic tent spaces in product domains must be defined.
We leave this interesting problems to jnterested readers.
The material of this survey related with the embedding theorems appeared in particular in the following way, many sharp new results embeddings were provided in recent papers of the first author and M.~Arsenovic, in context of harmonic function spaces of several complex variables. Then the first author extended completely these sharp embedding results for multyfunctional function spaces to analytic spaces and various other domains using similar type techiques and inthis survey we collected all such results. Various related new embeddings (sharp) but related with the trace function in various domains were also included in this survey.

In \cite{1} large amount of open problems were discussed in various  area Nevanlinna type function  spaces of several complex variables. We also refer our interested readers to \cite{1} for various interesting  issues related with embeddings in such type function spaces of several variables. Many interesting problems are open till now in this research area (see \cite{1} and also various references there).

Many old and new embedding theorems even in one dimension in the unit disk have great amount of applications in complex function theory and operator theory. We mention problems related with Toeplitz operators, Hankel operators, issues related with  description of continious linear  functionals and  various other problems. In higher dimension the situation is similar.
We hope this survey will serve as a valuable source of information for experts of function theory. In this survey we omit issues related with the embeddings in analytic function spaces  with
so called Bonami-Bekolle weights in one and higher dimension.
We refer the reader to papers mentioned in references concerning this new interesting issues.

We collected in the list of references various other papers of the first author and his coauthors. There some other new interesting sharp embedddings in various domains and analytic function  spaces of several complex variables or In harmonic function spaces can be seen. We refer the reader to these research papers.

We give in addition some very short comments concerning research papers on embeddings in analytic function spaces in one and several variables in various domains which can be seen in references of this paper and which we didnt mention in the text of this expository paper above (or we mentioned partially). In a paper of the first author in \cite{d39} various interesting embeddings  for so called holomorphic $F^{pq}_s$ spaces in the polydisk of spaces of analytic functions of Lisorkin-Triebel type can be seen. Similar type results in context of the unit ball and bounded stringly pseudoconvex domains can be seen in papers of J. Ortega and J. Fabrega \cite{d12}, \cite{d61}. In a paper of the first author with M.Radnia \cite{d23} new sharp interesting embeddings for some new analytic spaces in the unit ball can be seen. 

In a recent papers with Wen Xu \cite{d21}, \cite{d22} of the first author we can see new interesting embeddings for some analytic spaces based on Lorents spaces in the unit disk. We refer for similar type embedding results to a paper of the first author with A. Akbar (see \cite{d44}, \cite{d45}). In papers of D. Abate and coauthors \cite{d14}-\cite{d16} certain new interesting  embedding theorems in bounded strongly pseudoconvex domains can be seen. In older papers of J. Cima and D.~Luecking and coauthors (see \cite{d5}, \cite{d17}) the unit ball, unit disk cases were considered for certain Bergman type spaces.
In a recent research paper of the first author with N. Makhina \cite{d37} the first authors embedding results previously obtained with O. Mihic were extended to rather general domains in higher dimension for some mixed norm Bergman type analytic function spaces.

In a interesting paper of I. Verbitsky \cite{v1} a large amount of sharp embeddings for mixed norm analytic functoin spaces in the unit disk using delicate techniques were provided. In a series of recent papers of the first author with M. Arsenovic (see \cite{3}, \cite{d24}-\cite{k3}) harmonic function spaces of several variables were considered and many new interesting  embeddings theorems  there were given we partially indicated them above. In recent interesting research papers of the first author with R. Zhao, S. Kurilenko, O. Mihic and E. Povprits (see \cite{d332}, \cite{d39}, \cite{d29}-\cite{d331}) mixed norm spaces and various Bergman type and Herz type spaces were considered new embeddings were provided in context of unit ball and more difficult domains in higher dimension.

Multifunctional new analytic spaces also were considered in these research papers as well as in recent papers of the first author with E. Tomashevskaya (see \cite{d13}, \cite{d35}). Here new sharp interesting embeddings for some new multifunctional analytic function spaces in various complicated domains were found by the mentioned authors also.
In a recent interesting  paper of B. Sehba and D. Bekolle and various coauthors embeddings for analytic function spaces of Hardy and Bergman type were given in context of tubular domain over symmetric cone. These domains are unbounded and very general and many open interesting  problems here exits.
 
A recent paper in Filomat of the first author with M. Arsenovic \cite{d24} contains also such type results in tube domains over symmetric cones. We refer the reader also to the list of references of mentioned two last interesting research papers. In a paper of the first author with R. Zhao various embeddings were also provided.in the polydisk and unit ball. Finnaly, in older papers of F. Shamoyan (see \cite{d7}, \cite{d8}-\cite{d10}) various interesting embeddings can be seen in context of spaces of analytic functions the unit disk and polydisk.

Certain new interesting embeddings for new mixed norm analytic Lizorkin-Triebel spaces in the polydisk. Can be seen in \cite{d46}, \cite{d47}.

Papers on embeddings of first author with K. Avetisyan \cite{d48} can be viewed as an extention of the well known classical Littlewood-Paley inequality. Paper with R. Zhao \cite {d39} can be seen as embeddings in analytic spaces based on strong factorization of analytic function spaces in the unit disk.

In recent papers with O. Mihic and S. Maksakov of first author published in 2017 in ROMAI (see \cite{6}, \cite{11})  complete extentions of two sharp embedding theorems provided earlier by F. Shamoyan in the unit disk, polydisk  to more difficult tubular domains over symmetric conesand bounded strongly  pseudoconvex domains were given. The role of so called $r$-lattices in all theorems in all mentioned domains were very vital. These results probably are valid also in bounded symmetric domains with very similar proofs. We leave this for interested readers.

In a recent paper with Ali Akbar of first author analytic function spaces in the unit disk related with the so called p-Carleson measures were considered. The reader may find this paper in the list of references (see \cite{d45}). And new interesting  sharp embeddings were provided for such new function spaces.It will be interesting to extend these sharp embeddings to more complicated domains and  to analytic spaces of several complex variables. We pose this as a problem for interested readers.

In a recent short interesting  note with E. Ermakova \cite{d36} after carefull examination of proofs of some sharp embeddings provided earlier with S. Maksakov and O. Mihic in analytic spaces with product domains with mixed norm a general scheme was proposed y authors to get more general results of this type. Here many open interesting new problerms also exits. We leave this to interested readers. 

Many embeddings for various spaces in the ball can be seen in \cite{d49}. A natural problem here to extend these results to multifunctional case.

In $R^n$ many embeddings can be seen in \cite{d50}, \cite{d51}.

We refer the reader to \cite{d57}, \cite{d58} and references there for various recent results on analytic spaces embeddings in the unit ball and unit polydisk.

We refer the interested  readers  to papers of first author with S.Li \cite{k12}-\cite{k16} for various embeddings for so called myltifunctional analytic spaces in the unit ball and unit polydisk.

Finally, some interesting results on a topic of our expository paper can be seen in \cite{k1}-\cite{k11}.

We thank Natalia Makhina for her various consultations and her technical help during the preparation of this paper.

\end{document}